\documentclass[11pt]{article}
\usepackage{amscd,amsmath,amsxtra,amssymb,theorem,latexsym,amsfonts}
\usepackage[all]{xy}
\usepackage{epic,eepic}

\usepackage{oldgerm}
\usepackage{bbm}

\def\conv{\text{conv}}
\def\inter{\text{int}}

\newtheorem{satz}{Satz}
\newtheorem{lemma}[satz]{Lemma}
\newtheorem{theorem}[satz]{Theorem}

\newtheorem{corollary}[satz]{Corollary}
\newtheorem{prop}[satz]{Proposition}

\newcommand{\R}{{\mathbb R}}
\newcommand{\C}{{\mathcal C}}
\newcommand{\kf}{{\mathcal F}}
\newcommand{\ks}{{\mathcal S}}

\newcommand{\qed}{{\begin{flushright}$\Box$\end{flushright}}}

\CompileMatrices
\SelectTips{cm}{11}

\begin{document}

\author{Alexander Getmanenko}

\title{ \begin{center} \Huge Helly-type Theorems for Plane Convex Curves. 
\end{center} }

\maketitle

\begin{abstract} Families of translates and homothets of strictly convex curves are proven to possess Helly-type properties  generalizing  those of a circle. 
Weaker results are shown for arbitrary convex curves. 
\end{abstract}

This paper is concerned with proving a couple of Helly-type theorems for translates and homothets of convex curves on the Euclidean plane. It turns out that these Helly-type properties are very much like that of a circle and, as was amazing enough for me, admit of relatively simple and elementary proofs, which I 
present here.

Throughout this paper, a {\bf convex curve} $C$ is defined to be the boundary of a 
proper convex subset $B$ of Euclidean plane $\R^2$, i.e. $C=\partial B$.
A {\bf strictly convex} curve is a convex curve containing no line segment.

Given $C\subset \R^2$, we call any set of the form $\lambda C + v$, where $\lambda > 0$ and $v\in \R^2$ its {\bf homothet} and, specifically, its {\bf translate}, if $\lambda = 1$.

% notation : parallelogramm $abcd$

Classically (e.g., \cite{DGK}), the {\bf Helly number} of a family $\kf$ of sets is the minimal number $\chi$ satisfying the following property:\\
{\it if in a subfamily $F\subset \kf$ any subfamily $F_0\subset F$ with $\sharp F_0\le\chi$ has a nonempty intersection, then the whole $F$ has nonempty intersection. }

Helly theorem states that the Helly number of the family of all closed convex sets on the plane is 3. The Helly number of the family of all plane circles is 4, which is known from high-school geometry.

We define the {\bf translation order} $\mu_{Tr}$ of $C$ to be the minimal number $\mu$ such that: \\
{\it for any subset $S\subset \R^2$ such that for any $S_0\subset S$ with $\sharp S_0 \le 0$ points there is a translate $C'$ (depending on $S_0$) of $C$ containing $S_0$, then there exists a translate of $C$ containing the whole 
set $S$ }\\
Analogously the {\bf homothety order}  $\mu_{Ht}$ of $C$ is defined.

"Duals" of the facts above are that the translation order of every closed convex set on the plane is 3 and that of a circle is 4.

Our aim is to prove the following three theorems.

\begin{theorem} {\bf (i)} If $C$ is a convex plane curve, then 
$$ \mu_{Tr}(C) , \mu_{Ht}(C) \le 6.$$
{\bf (ii)} If $C$ is a strictly convex curve, then 
$$ \mu_{Tr}(C) , \mu_{Ht}(C) \le 4.$$
{\bf (iii)} If $C$ is strictly convex and bounded, then
$$ \mu_{Tr}(C) , \mu_{Ht}(C) = 4;$$
if $C$ is strictly convex and unbounded, then 
$$ \mu_{Tr}(C) = 3.$$
\end{theorem}

\begin{theorem} Let $C$ be a convex curve. Then:\\
{\bf (i)} the Helly number of the family of all translates of $C$ is $\le 6$;\\
{\bf (ii)} if $C$ is strictly convex, then this Helly number is $\le 4$; \\
{\bf (iii)} If $C$ is strictly convex and every 3 sets in a family of its translates $\C$ have a point in common and $\sharp\C \ge 5$ then the whole $\C$ has nonempty intersection. Moreover, the number 5 can be replaced by a smaller one for no bounded curve $C$.\\
{\bf (iv)} If $C$ is strictly convex and unbounded, 
then the Helly number of $\C$ is 3.
\end{theorem} 

\begin{theorem} Let $\C$ be a family of homothets of a strictly convex curve 
$C$. \\
{\bf (i)} If every 3 sets in $\C$ have a nonempty intersection, then the entire $\C$, except, maybe, one set, has a nonempty intersection.\\
{\bf (ii)} The Helly number of $\C$ is $\le 4$ and it is exactly $4$, if $C$ 
is bounded. 
\end{theorem}

All proofs are done simultaneously and are split into several statements, which are shown below.

The simple observation that follows will often be used in subsequent arguments.

\begin{lemma} \label{l1} Suppose for every set $S$ in a class $\ks$ of subsets of $\R^2$ there exists a subset $S_0$ with $\sharp S_0\le N$ such that there are $\le k$ translations (resp., homotheties) mapping $S_0$ into $C$. Then the translation (resp, homothety) order of $C$ with respect  
to the class $\ks$ is $\le k+N$.
\end{lemma}

The words {\it with respect to the class $\ks$} mean that we restrict ourselves to the sets of class $\ks$ in the definition of translation order. As you will see later, such classes can be, e.g., sets containing three collinear points. 
We subdivide the class of all subsets of a plane into several subclasses, and prove and estimate of a translation order for each of them. This idea is already used in the proof of the induction step in the proof of this lemma.

{\bf Proof.} Since the proof is much more general than the statement, let us do it for translates. 

Induction on $k$. Consider $S\in \ks$ such that any its $\le(N+k)$ -element subset can be mapped into $C$ by a translation.

If $k=1$, then take all $N$ points in $S_0$ and any point $p\in S\setminus S_0$.
Then there exists a translation $g_p$ such that $S_0\cup\{ p \} \subset g_pC$. But since $k=1$, the assumption of the lemma implies that $g_p=g_q$ for any $p,q\in S\setminus S_0$, and the whole $S$ is contained in $g_pC$.

Assume the statement be true for $(k-1)$, and let us show that it is true for $k$. Fix a subset $S_0$ and $k$ translates $C_1,\dots,C_k$ of $C$ containing $S_0$. If $S\subset \bigcap_i C_i$ for some $i$, we are done, otherwise take
a point $p\in S\setminus X_i$. Then there are at most $(k-1)$ translates of $C$ containing the $(N+1)$-elements subset $S_1=S_0\cup\{p\}$, and we may apply the induction hypothesis to the class $\ks_1=\{ S\in \ks \ : \ S_0\not\subset  \bigcap_i C_i \}$. \qed

\begin{lemma} \label{collin} {\bf (i)} If $C$ contains points $x,y,z,x+v,y+v,z+v$, then two of
the four
vectors $(x-y),$ \linebreak $(x-z),(y-z),v$ are collinear. \\
{\bf (ii)} Suppose $0$ is the origin, $\lambda\ne 1$ is a positive number, and $C$ contains six points $x,y,z,\lambda x, \lambda y, \lambda z$. Then there are tree collinear points among $x,y,z,0$. 
\end{lemma} 

{\bf Proof.} {\bf (i)} Let $\pi$ be a projection map in the direction of the vector $v$. We may assume that $\pi(x),\pi(y),\pi(z)$ are pairwise distinct,  otherwise the lemma holds trivially. Without loss of generality, $\pi(y)$ lies between $\pi(x)$ and $\pi(z)$. If $x,y,z$ are not collinear, then either $y$ or $y+v$ lies inside the parallelogramm $xy(x+v)(y+v)$, contradicting to $C$ being a convex curve.

{\bf (ii)} Suppose $x,y,z$ are not collinear. Then $0 \not\in \conv(x,y,z)$, i.e.
$0$ and $x,y,z$ are separated by a line $\ell$. Let $\pi$ denote the projection
from the center $0$ to the line $\ell$. If any two points among $\pi(x),\pi(y),\pi(z)$ coincide, the statement is fulfilled; if not, let $\pi(y)$ lie between $\pi(x)$ and $\pi(z)$. As $x,y,z$ are not collinear, either $y$ or $\lambda y$ lies inside trapezoid $xy(\lambda x)(\lambda y)$ that contradicts to convexity of the curve $C$. \qed

\begin{prop} \label{Mule6} $\mu (Tr,C), \mu(Ht,C) \le 6$.
\end{prop}

{\bf Proof.} We do the proof for homothets; the proof for translates is analogous.

Let $S\subset\R^2$, $\sharp S \ge 6$ be a set such that every subset of $\le 6$ points can be moved inside $C$ by a positive homothety. Suppose all points of $S$, except, maybe, a point $p$, belong to a line $\ell$. Choose a subset $S_0$
consisting of 4 or 5 points in S: the leftmost and the rightmost points of the intersection $S\cap \ell$, $p$ (if it exists) and  any two other points of $S$. Let $gC$ be a homothet of $C$ containing $S_0$, then the entire set $S$ is  contained in $gC$.

Now assume that $S$ contains at least 4 points in general position, say, $a,b,c,d$ being vertices of a convex 4-gon. 

\begin{lemma} If the points $a,b,c,d$ are in general position and define a convex 4-gon, then all the homothets of $abcd$ that can be embedded into $C$  have the same homothety center.
\end{lemma}

{\bf Proof.}
Suppose $abcd$ and $g(abcd)=a'b'c'd'$ are in $C$ for some homothety $g$, and $0$ is the homothety center of $g$. If $g$ is a translation, formally assume $0$ to be a point at infinity, and we still may use lemma \ref{collin} part (ii).

Since $abc$ is homothetic to $a'b'c'$, and $a,b,c$ are not collinear, we have that one of the triples $(a,b,0)$, $(a,c,0)$ or $(b,c,0)$ is collinear (by lemma \ref{collin}) and, correspondingly, one of the 4-tuples   
$(a,b,a',b')$ or $(a,c,a',c')$ or $(b,c,b',c')$ is collinear. Analogously, if $abd$ is homotetic to $a'b'd'$, we get collinearity of one of $(a,b,a',b')$, $(a,d,a',d')$, $(b,d,b',d')$.

If $a,b,a',b'$ are collinear, then, as $acd$ is homothetic to $a'c'd'$, one of the 4-tuples  ${a,c,a',d'}$, ${a,d,a',d'}$, ${c,d,c',d'}$ is collinear. In fact, the first two 4-tuples cannot be collinear, because then there are two lines containing the point $a$ and invariant under the homothety $g$, hence $0=a$; but then one of the points $c$ or $c'$ lies inside the convex hull of the 8 points $a,b,c,d,a',b',c',d'$, contradiction. Hence ${c,d,c',d'}$ is collinear  and $0=(ab)\cap (cd)$. 

Suppose $a,b,a',b'$ are not collinear, and ${a,c,a',c'}$ are collinear. Then, doing the argument analogous to the one above (with $abd$ instead of $acd$), we see that $0=(ac)\cap (bd)$, i.e. $0 \in \inter abcd$, and this is impossible.
 
If neither ${a,b,a',b'}$ nor ${a,c,a',c'}$ is collinear, then ${b,c,b',c'}$ is collinear, and one of the 4-tuples ${a,d,a',d'}$ or ${b,d,b',d'}$ is collinear. Analogously one can see that the latter cannot be collinear, and $0=(bc)\cap (ad)$.

We have proven that either $0=0_1=(bc)\cap (ad)$ or $0=0_2=(ab)\cap (cd)$. Suppose now $a_1b_1c_1d_1$ and $a_2b_2c_2d_2$ are homothetic to $abcd$ with centers $0_1$ and $0_2$ respectively. A reader is urged to draw a picture to convince oneself, then, that a convex curve $C$ cannot contain all the 12 points
 $a,b,\dots,d_2$.
This contradiction proves the lemma. \qed

We are now  resuming the proof of proposition \ref{Mule6}
Suppose $abcd$ may be homothetically embedded into $C$ uniquely, then by lemma 
\ref{l1} the result follows immediately. Otherwise let $0=(bc)\cap (ad)$ be the 
common homothety center of all its homothetic images. Then the lines $\ell=(bc)$ and $k=(bc)$ are supporting lines of $C$. If $S\subset k\cup \ell$, then let us 
choose 3 points on either of the lines, including the rightmost and the leftmost points of $S\cap \ell$ and of $S\cap k$, respectively. If there is a point 
$x\in S\setminus (k\cup\ell)$, then there exists only one homothet of $C$ containing 
$a,b,c,d,x$, and by lemma \ref{l1} the result follows. \qed

Let us briefly discuss another way of proving this lemma for translates. We will study the number of translates of a nondegenerate triangle $T={x,y,z}\subset S$ contained in $C$. Let $V\subset \R^2$ be the set of all $v$ such that $T+v\subset C$. Clearly $0\in V$. Assume V contains two linearly independent vectors $u,v$. By lemma \ref{collin} we have after, possibly, changing the notation for $x,y,z$, that $u \parallel (x-y)$, $v \parallel (x-z)$ and $(u-v) \parallel (z-y)$. This means that the triangles $xyz$ and $0uv$ are homothetic. Then the lines connecting the points $x,y+u,z+v$ are supporting lines to $C$, and $C$ is the triangle with the  vertices at $x,y+u, z+v$, and $C$ contains precisely 3 translates of $C$. If all vectors in $V$ are collinear, proceed as in the proof of the proposition.

The proof of the proposition \ref{Mule6} and lemma \ref{collin} yield:

\begin{prop} \label{Mule4} If $C$ is a strictly convex curve, then there exists no more than one homothet (in particular, at most one translate) of $C$ containing a given three points. Therefore,
$$ \mu_{Tr}(C), \mu_{Ht}(C) \ \le 4. $$  
\end{prop} 

\begin{corollary} Let $C$ be a convex curve. If any $m=6$ sets in a family $\C$ of translates of $C$ have a nonempty intersection, then the whole family $\C$
has nonempty intersection. Moreover, if $C$ is strictly convex, we may assume $m=4$.
\end{corollary}

{\bf Proof} Let $V$ be defined by $\C=\{ x+C | x\in V\}$. A subfamily 
$\{ x+C | x\in V_0 \} \subset V$ has a nonempty intersection iff for some $z$ $V_0\subset z-C$. The assumption implies that any $m$ points of $V$ are contained in a translate of $(-C)$, hence, the whole $V$ contains in a translate of $(-C)$, that proves the corollary. \qed

\begin{corollary} Any two distinct positive homothets of a strictly convex curve $C$ have at most two intersection points.
\end{corollary} 

{\bf Proof.} If we had $\ge 3$ intersection points, they would define a triangle that can be embedded into $C$ in two different ways, which contradicts the proposition \ref{Mule4}. \qed

\begin{corollary} Let $C$ be a strictly convex curve. \\
(i) if every $m=3$ sets from a family $\C=\{C_i\}$ of positive homothets of $C$ have a nonempty intersection,then the whole family $\C$, except, perhaps, one set, have a nonempty intersection.\\
(ii) if any $m=4$ sets from a family $\C=\{C_i\}$ of positive homothets of $C$ have nonempty intersection, then the whole family $\C$ has nonempty intersection.
\end{corollary}

{\bf Proof} Suppose any 3 sets in $\C$ have an intersection point. If there are two sets of $\C$ that intersect in only one point, then this point belong to all other curves $C_i$ as well. Therefore we may assume that any two sets in $\C$ intersect in exactly two points. 

Denote by $\C_0$ a maximal (with respect to inclusion) subfamily of $\C$ with a common point, say, $p$. Clearly $\sharp \C_0|\ge 3$. Let $T=\{ q\ne p \ : \ q\in C_i\cap C_j , C_i,C_j\in \C_n\}$. Suppose $D\in \C \setminus \C_0$, then $D\supset T$. 

If $\sharp T\ge 3$, then there can be no more than one curve $C'$ containing $T$. 

If $T=\{ x \}$, then each set in $\C_0$ contains $x$, hence $\sharp\C\le 2$, which contradicts to $\sharp\C\ge 3$. 

If $T=\{ x,y \}$, then $\C_0=\C_x\cup \C_y$, where $\C_x=\{ C_i\in \C_0 \ : \ C_i\ni x  \}$, $\C_y=\{ C_i\in\C \  : \ C_i\ni y \}$. By definition of $T$ we get $\sharp C_x, \sharp C_y \ge 2$ and $\sharp(\C_x\cap \C_y)\le 1$ because any homothet of $C$ in $\C_x\cap\C_y$ must contain 3 distinct points $p,x,y$. Therefore there exists $D_x\in \C_x\setminus\C_y$, $D_y\in \C_y\setminus \C_x$. We have $D_x\cup D_y=\{p,z\}$, with $z\in\{ x,y\}$, say, $x=z$, $D_y\in\C_x$ contradicting to the choice of $D_y$. Therefore the case $\sharp T = 2$ is impossible, and the proof of (i) is complete. 

(ii) By the argument above we may assume $\sharp T\ge 3$, $\{ x,y,z \}\subset T$. Choose 3 (not necesserily distinct) sets in $\C_0$: $E_x\not\ni x$, $E_y\not\ni y$, $E_z\not\ni z$. As $E_x\cap E_y\cap E_z = \{ p \}$, then $E_x\cap E_y\cap E_z \cap D = \{ p\}$, i.e. $p\in D$ and $\C=\C_0$. This proves part (ii). 
\qed

\begin{corollary} Let $C$ be a strictly convex curve. If any 3 sets in the family $\C$ of translates of $C$ have a point in common and $\sharp\C\ge 5$, then the entire family $\C$ has nonempty intersection.
\end{corollary}

{\bf Proof} The proposition and the corollary above reduce the question to the following: given 5 points $x_i$, $i=1,\dots,5$ such that every 3 of them are contained in a translate of $C$, then all the five points are contained in a translate of $C$.

Lemma \ref{collin} implies that there are at most two translates of a given pair of points $x_1x_2$ inside $C$. For each of the triple $x_1x_2x_i$, $i=3,4,5$ there exists precisely one translation of $\{ x_1x_2x_i\}$ putting a triple onto the curve $C$. Hence there is a translation that puts 4 point, say, $x_1,x_2,x_3,x_4$ onto $C$. Doing the same argument starting from the pair $x_1,x_5$, we see that (up to renumbering) $C$ contains a translate of the 4-tuple $x_1,x_2,x_3,x_5$. As the images of $x_1,x_2,x_3$ under these translations have to coincide, we get that this translation moves all the 5 points onto $C$. \qed

\begin{prop} For any compact convex curve there are 4 points belonging to no translate of $C$ such that every 3 of them lie on some translate of $B$.
\end{prop}

{\bf Proof,} Inscribe an affinely regular hexagon $x_1x_2x_3x_4x_5x_6$ with center $c$ in  $C$ and take points $x_1,x_3,x_5,c$ that do the job. \qed 

\begin{prop} The translation order of any stricly convex unbounded curve equals 3. \\ Equivalently, if $C$ is an unbounded stricly convex curve and $\C$ is a family of its translates, such that any 3 sets in $\C$ have a point in common, then all the family $C$ has nonempty intersection, and the number 3 cannot be replaced by a smaller number.
\end{prop}

{\bf Proof.} We proof the first formulation. Take any 3 collinear points $x_1,x_2,x_3$, then there obviously exists three translates of $C$ containing 2 of these 3 points each and having empty intersection. The rest follows by lemma \ref{l1} because there is at most one translation that puts 2 given points onto $C$. \qed   

Collecting the statements, we get the three theorem announced at the very beginning.

The results presented in this paper form a part of my Bachelor thesis defended at the Yaroslval State University in July, 1997. I thank my supervisor V.L.Dol'nikov who helped me a lot during my study in Yaroslavl and encouraged to prepare this publication.

\end{document}